\documentclass[12pt]{amsart}
\usepackage{amssymb}

\newtheorem{theorem}{Theorem}[section]
\newtheorem{lemma}[theorem]{Lemma}
\newtheorem{proposition}[theorem]{Proposition}
\newtheorem{corollary}[theorem]{Corollary}
\newtheorem{definition}[theorem]{Definition}

\theoremstyle{definition}
\newtheorem{remark}[theorem]{Remark}
\newtheorem{example}[theorem]{Example}
\numberwithin{equation}{section}

\newcommand{\Z}{\mathbb{Z}}
\newcommand{\R}{\mathbb{R}}
\newcommand{\C}{\mathbb{C}}

\newcommand{\CP}{\mathbb{CP}}

\newcommand{\F}{\mathcal{F}}
\newcommand{\Map}{\mathrm{Map}}
\newcommand{\Diff}{\mathrm{Diff}}
\newcommand{\mm}[1]{\begin{pmatrix}#1\end{pmatrix}}

\title[Factorizations in $SL(2,\Z)$ and inequivalent Stein fillings]{Factorizations
in $SL(2,\Z)$ and simple examples of inequivalent Stein fillings}
\author{Denis Auroux}
\address{Department of Mathematics, UC Berkeley, Berkeley CA 94720-3840, USA}
\email{auroux@math.berkeley.edu}
\thanks{The author was partially supported by NSF grants DMS-1007177 and
DMS-1264662.}

\begin{document}

\begin{abstract}
We give simple examples of elements of $SL(2,\Z)$ admitting
inequivalent factorizations into products of Dehn twists. This can be
interpreted in terms of inequivalent Stein fillings of the same 
contact 3-manifold by genus 1 Lefschetz fibrations over the disk.
\end{abstract}

\maketitle

\section{Introduction}

Lefschetz fibrations have risen to prominence in recent years as a
convenient way to describe symplectic 4-manifolds. In particular,
Lefschetz fibrations over the disk correspond to Stein
fillings of contact 3-manifolds; see e.g.\ \cite{AkOz,LP,OS03,SeBook}.

The classification of Lefschetz fibrations remains poorly understood to date, with a wealth
of ``exotic'' examples constructed in recent years. For instance, 
genus~2 (or higher) Lefschetz fibrations over the disk have been 
used to find contact 3-manifolds which admit infinitely many 
inequivalent Stein fillings; see e.g.\ \cite{AEMS,OS03}.
By contrast, the classification in genus~1
has generally been thought to be much simpler, perhaps due to
the classical result of Moishezon of Livne \cite{Moishezon} according to
which genus 1 Lefschetz fibrations over $S^2$ are holomorphic and
classified by their number of singular fibers.

In this paper, we show that genus 1 Lefschetz fibrations over the disk
are much more subtle than their closed counterparts. Specifically, we describe some 
simple examples of such fibrations which give different Stein fillings 
(e.g., with different first homology groups) of the same contact 
3-manifold. These arise from inequivalent factorizations of the same 
element in $SL(2,\Z)$ as a product of Dehn twists. These also lead to
various other interesting small examples, e.g.\ of
different symplectic submanifolds in $B^4$ filling the same 3-strand braid, or 
different Lagrangian disks in a Stein manifold bounding the same
Legendrian knot.

Our first and main example, with four singular fibers, is a pair of Lefschetz fibrations
that have already been studied in the context of mirror symmetry, where they
occur as the mirrors of Hirzebruch surfaces ($\CP^1$-bundles over $\CP^1$): 
i.e., for instance, the derived Fukaya 
categories of vanishing cycles of these Lefschetz fibrations are equivalent 
to the derived categories of coherent sheaves of the latter spaces \cite{AKO}.
The two fillings are distinguished by their first homology 
groups; see Proposition \ref{prop:exampleMS}.

These examples sit inside rational elliptic fibrations with $I_8$ singular
fibers (namely, as the complement of the $I_8$ fiber and of a section); the
existence of two distinct types of such rational elliptic fibrations is
well-known in algebraic geometry (as a consequence of the classification
of extremal fibrations \cite{Beauville,MP}, see e.g.\ \S\,VIII of \cite{M}). One can similarly look
at other examples of extremal or nearly-extremal elliptic fibrations
\cite{MP,MP2}, 
such as
elliptically fibered K3 surfaces with $I_{18}$ or $I_{19}$ singular fibers, or $E(3)$
elliptic surfaces with $I_{29}$ singular fibers. While these do give rise to
other examples, all those we found can be understood in 
terms of smaller building blocks, and the relevance of extremal elliptic
fibrations to the question at hand is far from clear.

The smallest possible examples one could hope for are genus 1 Lefschetz
fibrations with only three singular fibers (or even two, if one does not
require the fillings to be topologically distinct).
It turns out that such examples abound.
For instance, the example described in \S \ref{s:exampleMS} can be simplified
by discarding one of the singular fibers, at the 
expense of making its conceptual significance less clear. 
However, there exist many other
examples of inequivalent genus 1 Lefschetz fibrations with
three singular fibers and the same boundary monodromy; we list some of them
(found by a computer search) in \S \ref{s:example3fibers}.
Some of these examples can be distinguished by their homology. Others require
a more subtle invariant of Lefschetz fibrations that we describe in 
\S \ref{s:invt3fibers}.

Finally, we point out that various classification results can still be hoped
for in spite of these fairly discouraging examples. In genus 1, the mapping
class group has a fairly simple structure, and one can enumerate the
possible factorizations of a given element into a given number of Dehn
twists \cite{Cadavid}. In fact, when there are only three singular fibers
(and still in genus 1), the invariant described in \S \ref{s:invt3fibers} 
seems to capture nearly all the information.  In a different vein, it follows 
from the results in \cite{stabslf} that any two Stein
fillings of a given contact 3-manifold 
with the same Euler characteristic and signature
become equivalent under stabilization by performing the same sequence of
handle attachments at the contact boundary; see Theorem \ref{thm:stab}. Thus
the phenomena we discuss below are inherently ``unstable''.

\subsection*{Acknowledgements}
I would like to thank Carlos Cadavid for sending me his preprint
\cite{Cadavid}, which prompted me to think about this question, as well
as Paul \hbox{Seidel} for encouraging comments and for pointing out several 
interesting features of our main
examples, and the referees for their useful suggestions. This work was partially supported by NSF grants DMS-1007177 and
DMS-1264662.

\section{Lefschetz fibrations and monodromy factorizations}

\subsection{Lefschetz fibrations}
A Lefschetz fibration over the disk is a map $f:M^4\to D^2$ whose
smooth fibers are oriented surfaces, and whose only critical points
occur over the interior and are modelled on the complex Morse singularity
$(z_1,z_2)\mapsto z_1^2+z_2^2$ in an orientation-preserving coordinate
chart. The singular fibers of $f$ thus are obtained from a smooth fiber
$F$ by collapsing a simple closed curve, called the {\em vanishing cycle}, to an
ordinary double point; the monodromy
around each singular fiber is given by a right-handed Dehn twist about the
appropriate vanishing cycle.

The relation to symplectic geometry is the following. Assume the fiber $F$
has non-empty boundary, and is equipped with an exact symplectic
structure, in such a way that all the vanishing cycles are exact
Lagrangian submanifolds (this can always be arranged when the
vanishing cycles are nonzero in homology). The total space $M$ of the Lefschetz
fibration $f$ then carries an
exact symplectic structure, canonical up to deformation, while the 
restriction of $f$ to the boundary of $M$ (a contact 3-manifold) endows 
$\partial M$ with an open book decomposition which supports the contact structure.
Topologically, $M$ is obtained from $F\times D^2$ by attaching standard
Weinstein 2-handles along the vanishing cycles in parallel copies of the 
fiber in
$\partial(F\times D^2)$. See e.g.\ \cite{AkOz,SeBook} for more details.

Taking the reference fiber $F$ to lie over a base point near the boundary of 
$D^2$, and choosing a distinguished collection of paths that connect the base point to the 
various critical values of $f$ (assumed to be distinct), 
we obtain a {\em distinguished basis of 
vanishing cycles} $(\gamma_1,\dots,\gamma_r)$ in $F$.
The monodromies around the various singular fibers, i.e.\ the Dehn twists
$\tau_1,\dots,\tau_r$ about $\gamma_1,\dots,\gamma_r$, completely determine
the topology of the Lefschetz fibration $f$; moreover, their product $\phi$ is
the monodromy of the open book induced by $f$ on~$\partial M$. Thus, we
can describe $f$ by its {\em monodromy factorization},
i.e.\ a decomposition of $\phi$ into a product of Dehn twists
$\phi=\tau_1\cdot\ldots\cdot \tau_r$, in the mapping class group
$\Map(F,\partial F)=\pi_0\Diff^+(F,\partial F)$.

The braid group $B_r$ acts simply transitively on the set of distinguished
bases of paths; the corresponding action on monodromy factorizations is
called {\em Hurwitz equivalence}, and is generated by the Hurwitz moves
$$(\tau_1,\dots,\tau_i,\tau_{i+1},\dots,\tau_r) \sim (\tau_1,\dots,
\tau_i\tau_{i+1}\tau_i^{-1},\tau_i,\dots,\tau_r)\quad \text{for }
1\le i<r$$
and their inverses. In terms of vanishing cycles, this amounts to
replacing $\gamma_i$ and $\gamma_{i+1}$ by $\tau_i(\gamma_{i+1})$ and
$\gamma_i$ respectively. Hurwitz equivalence classes of monodromy
factorizations correspond to isomorphism classes of Lefschetz fibrations
with a marked fiber, i.e.\ with a fixed identification of $F$ with 
some abstract oriented surface with boundary. Changing this identification
by an element $\psi$ of the mapping class group amounts to replacing each
vanishing cycle $\gamma_i$ by its image $\psi(\gamma_i)$, i.e.\ to a
{\em global conjugation} of the monodromy factorization, replacing
each $\tau_i$ by $\psi\tau_i\psi^{-1}$. (Of course, this now yields a
factorization of $\psi\phi\psi^{-1}$.)
The classification of Lefschetz fibrations over the disk thus amounts to
that of monodromy factorizations in the mapping class group up to Hurwitz
equivalence and global conjugation (cf.\ e.g.\ \cite{AurouxMCG,stabslf}).

\subsection{The genus one case}
In this paper we will focus specifically on the case where $F$ is a torus with one boundary
component. The mapping class group of $T^2$ is $SL(2,\Z)$, while that
of a genus 1 surface with one boundary component is
$$\Gamma:=\Map_{1,1}=\widetilde{SL}(2,\Z),$$ a central
extension of $SL(2,\Z)$ by $\Z$ which can be represented as the preimage of
$SL(2,\Z)$ in the universal cover of $SL(2,\R)$ (hence the notation). Because every
punctured elliptic curve is a double cover of
the complex plane branched at three points, $\Gamma$ is also isomorphic to
the 3-strand braid group $B_3$. 

The group $\Gamma$ is generated by the Dehn twists $a$ and $b$
about simple closed curves $\alpha,\beta$ representing the two
$S^1$ factors of the torus, with the relation $aba=bab$.
The {\em boundary twist}
$\delta=(ab)^6$ (i.e., the Dehn twist about a
boundary-parallel curve) is central and generates the kernel of the quotient map $\Gamma\to SL(2,\Z)$.
Since Dehn twists in $\Gamma$ are determined by their images in $SL(2,\Z)$,
a monodromy factorization in $\Gamma$ is specified unambiguously by its image in
$SL(2,\Z)$, a fact that we will use repeatedly in the next sections.
Moreover, two factorizations of the same element of $SL(2,\Z)$ into products
of the same numbers of Dehn twists lift to factorizations of the same
element in $\Gamma$. 
(Both properties follow from the observation that
an element of $\Gamma$ is determined by its images in $SL(2,\Z)$ and in
the abelianization $\mathrm{Ab}(\Gamma)\simeq \Z$; under the latter map,
Dehn twists map to 1 while the central element $\delta$ maps to 12).

To be more explicit, the generating Dehn twists $a$ and $b$ map to the two generators
$$A=\mm{1&1\\0&1}\quad \text{and}\quad B=\mm{1&0\\-1&1}$$
of $SL(2,\Z)$; more generally, the Dehn twist $\tau_{p,q}$ about a simple closed curve
representing the class $p[\alpha]+q[\beta]=(p,q)\in H_1(F)\simeq\Z^2$ corresponds to
the matrix
$$T_{p,q}=\mm{1-pq & p^2\\ -q^2 & 1+pq}\in SL(2,\Z).$$

\section{The main example}\label{s:exampleMS}

\begin{proposition}\label{prop:exampleMS}
The monodromy factorizations 
\begin{equation}\label{eq:mir_F1}
\phi=\tau_{-3,1}\cdot\tau_{0,1}\cdot\tau_{3,1}\cdot\tau_{1,0}=
(a^{-3}ba^3)\cdot b\cdot (a^3ba^{-3})\cdot a \quad \text{and}
\end{equation}
\begin{equation}\label{eq:mir_F0}
\phi=\tau_{-2,1}\cdot\tau_{0,1}\cdot\tau_{0,1}\cdot\tau_{2,1}=
(a^{-2}ba^2)\cdot b\cdot b\cdot (a^2ba^{-2})
\end{equation}
of $\phi=a^{-8}\delta$ in $\Gamma$ define inequivalent
genus 1 Lefschetz fibrations $f_1,f_2$ over the disk. The corresponding
Stein fillings $M_1,M_2$ of the open book with monodromy $\phi$ are
distinguished by their first homology groups: $H_1(M_1,\Z)=0$ while
$H_1(M_2,\Z)=\Z/2$.
\end{proposition}

\proof
The identities \eqref{eq:mir_F1} and \eqref{eq:mir_F0} can be checked
either by direct calculation in $\Gamma=\langle a,b\,|\,aba=bab\rangle$,
or by working in $SL(2,\Z)$, where it is easy to verify that
\begin{eqnarray*}
T_{-3,1}T_{0,1}T_{3,1}T_{1,0}&=&\mm{4&9\\-1&-2}\mm{1&0\\-1&1}\mm{-2&9\\
-1&4}\mm{1&1\\0&1}=\mm{1&-8\\0&1}\\
\mathrm{and}\quad T_{-2,1}T_{0,1}T_{0,1}T_{2,1}&=&
\mm{3&4\\-1&-1}\mm{1&0\\-1&1}^2\mm{-1&4\\-1&3}=\mm{1&-8\\0&1}.
\end{eqnarray*}
The Lefschetz fibrations $f_1$ and $f_2$ are easily distinguished by the
fact that the vanishing cycles of $f_1$ generate $H_1(F)\simeq \Z^2$
while those of $f_2$ only generate an index 2 subgroup. Accordingly, the
first homology groups of $M_1$ and $M_2$, which are isomorphic to the
quotients of $H_1(F)$ by the span of the vanishing cycles, are also
different.
\endproof

Another way to distinguish the monodromy factorizations \eqref{eq:mir_F1}
and \eqref{eq:mir_F0} in $\Gamma\simeq B_3$ is to consider their images in 
$SL(2,\Z/2)\simeq S_3$: while the factors in \eqref{eq:mir_F1} generate the
whole group, those in \eqref{eq:mir_F0} all map to the same element.

\begin{remark}
Viewed as a pair of factorizations in the braid group $B_3$, this example
can be thought of as a simpler analogue of that given in
\S 5 of \cite{AKS} (which involves products of 6 half-twists in
$B_4$). In fancier language, the 
``generalized Garside problem'' (i.e., whether a factorization is determined
by the product of its factors) also has a negative answer for products of
half-twists in $B_3$. 

Viewed as braid group factorizations, \eqref{eq:mir_F1}
and \eqref{eq:mir_F0} determine properly embedded smooth
symplectic surfaces $\Sigma_1$ and $\Sigma_2$ in the 4-ball, whose
boundary is the same transverse link. Namely, $\Sigma_1$ and $\Sigma_2$
(which can in fact be chosen algebraic)
are characterized up to isotopy by the requirement that projection to the
first two coordinates makes $\Sigma_i$ a 3-sheeted
branched cover of the disk, with four simple branch points around which
the monodromies are given by the factors in \eqref{eq:mir_F1} resp.\
\eqref{eq:mir_F0}; see e.g.\ \cite{Orevkov,KK,AKS}.
The symplectic surfaces $\Sigma_1$ and $\Sigma_2$ are easily distinguished 
by the fact that $\Sigma_1$ is connected while $\Sigma_2$ is not. 
In this language, the symplectic 4-manifolds $M_1$ and $M_2$ are the 
double covers of $B^4$ branched at $\Sigma_1$ and $\Sigma_2$ respectively.
We note that the trick used by Geng \cite{Geng} to modify the example of
\cite{AKS} into a pair of connected symplectic surfaces distinguished
by the fundamental groups of their complements fails in this example,
as $\pi_1(B^4\setminus \Sigma_2)\simeq \Z^2$ is too small (namely, 
the fundamental groups of the complements would be quotients of $\Z^2$
hence abelian, but for a smooth connected surface the first homology group 
of the complement is always $\Z$).
\end{remark}

The Lefschetz fibrations $f_1$ and $f_2$
are closely related to the toric Landau-Ginzburg mirrors of the Hirzebruch
surfaces $\mathbb{F}_1$ ($\CP^2$ blown up at one point) and
$\mathbb{F}_0=S^2\times S^2$ (or equivalently up to deformation,
$\mathbb{F}_2=\mathbb{P}(\mathcal{O}_{\mathbb{P}^1}(-2)\oplus
\mathcal{O}_{\mathbb{P}^1})$) respectively; see \S 5~of~\cite{AKO} for a discussion of
these examples from the perspective of homological mirror symmetry. 
Specifically, the toric mirror of $\mathbb{F}_k$
($k=0,1,2$),
given by the Laurent polynomial $W_k=x+y+x^{-1}+x^{-k}y^{-1}:(\C^*)^2\to \C$, 
is an elliptic fibration over the complex plane, whose
fibers have four punctures instead of one. Modifying the fibrations $W_1$
and $W_2$ by partial 
fiberwise compactification (to have once-punctured tori as fibers) and choosing all
the vanishing cycles to be exact, we obtain mirrors of $\mathbb{F}_1$ and
$\mathbb{F}_2$ which are exactly the Lefschetz fibrations $f_1$ and $f_2$.

From another perspective,
there are two different types of rational elliptic fibrations with one
$I_8$ singular fiber and four ordinary ($I_1$) nodal singular fibers
\cite{MP,M};
$f_1$ and $f_2$ can be constructed from these by deleting the $I_8$ fiber
and a section. 

The difference between $f_1$ and $f_2$ disappears after adding just one
new singular fiber with the same vanishing cycle $\alpha$ to both 
of them (i.e., adding an extra factor $\tau_{1,0}=a$ to both \eqref{eq:mir_F1} and \eqref{eq:mir_F0}):

\begin{lemma}
The factorizations $a^{-7}\delta=a\cdot (a^{-3}ba^3)\cdot b\cdot
(a^3ba^{-3})\cdot a$ and 
$a^{-7}\delta=a\cdot (a^{-2}ba^2)\cdot b\cdot b\cdot (a^2ba^{-2})$
are Hurwitz equivalent.
\end{lemma}

\proof We perform successive Hurwitz moves on the first expression,
moving the underlined factors across their neighbors (which undergo
conjugation) each time:
\begin{eqnarray*}
\underline{a}\cdot (a^{-3}ba^3)\cdot b\cdot (a^3ba^{-3})\cdot \underline{a} 
&\sim& (a^{-2}ba^2)\cdot \underline{a}\cdot \underline{b}\cdot a\cdot (a^2ba^{-2})\\
&\sim& (\underline{a^{-2}ba^2})\cdot \underline{b}\cdot a\cdot b \cdot (a^2ba^{-2})\\
&\sim& a\cdot (a^{-2}ba^2)\cdot b\cdot b\cdot (a^2ba^{-2}).
\end{eqnarray*}
\vskip-3ex \endproof

Thus, attaching Weinstein 2-handles to $M_1$ and $M_2$ along the
same Legendrian knot in the boundary (note that $\partial M_1=\partial M_2$ 
as contact manifolds) yields
new Stein manifolds $M_i^+=M_i\cup_\partial H_i$ (carrying Lefschetz
fibrations $f_i^+:M_i^+\to D^2$ with five singular fibers) which are 
deformation equivalent: 
$M_1^+\simeq M_2^+$, and we denote this manifold simply by $M^+$.  

As pointed out by Paul Seidel, this implies:

\begin{corollary}
There exists a Legendrian knot $K\subset \partial M^+$ which admits two
non-isotopic fillings by properly embedded Lagrangian disks 
$D_1,D_2\subset M^+$, $\partial D_i=K$. The two fillings $D_1,D_2$ are distinguished
by the first homology group of their complements.
\end{corollary}

\proof
Take $D_i$ to be the co-core of the Weinstein handle
$H_i$ in $M_i^+$, or in other terms, the Lefschetz thimble
associated to a vanishing path that runs from the critical point of
$f_i^+$ which lies inside $H_i$ straight to
a base point $q\in \partial D^2$.  The boundaries of $D_1$ 
and $D_2$  are the same Legendrian knot in $\partial M_1^+=\partial M_2^+$, 
namely the loop $\alpha$ inside the fiber 
$(f_1^+)^{-1}(q)=(f_2^+)^{-1}(q)$. Indeed, since the monodromy
factorizations of the Lefschetz fibrations $f_1^+$ and $f_2^+$ are 
Hurwitz equivalent, the isomorphism between them is compatible with the
chosen markings of the reference fibers, and in particular maps $\alpha$ in
$(f_1^+)^{-1}(q)$ to the corresponding loop in $(f_2^+)^{-1}(q)$.
On the other hand, since
$D_i$ is the co-core of the handle $H_i$, its complement $M_i^+\setminus D_i$ retracts
onto $M_i$, and so $H_1(M_1^+\setminus D_1)=0$ while
$H_1(M_2^+\setminus D_2)\simeq \Z/2$, by
Proposition \ref{prop:exampleMS}.
\endproof

\begin{remark}
The fact that $f_1$ and $f_2$ become isomorphic after adding a new singular 
fiber with monodromy $a$ to each of them corresponds under mirror symmetry 
to the classical fact that blowing up a point on $S^2\times S^2$ yields the
same del Pezzo surface as blowing up two points on $\CP^2$. Namely, as
shown in \cite{AKO2}, blowing up a del Pezzo surface modifies its
mirror by adding an extra vanishing cycle which is
``pulled from the fiber at infinity''; in our case these vanishing cycles
represent the class $[\alpha]$, and passing from $f_k$ to $f_k^+$
amounts to passing from the mirror of a Hirzebruch surface 
to that of its blowup.
\end{remark}

\section{Examples with two or three singular fibers}\label{s:example3fibers}

Modifying \eqref{eq:mir_F1} by a single Hurwitz move, we can rewrite it as
$(a^{-3}ba^3)\cdot b\cdot a\cdot (a^2ba^{-2})$, which makes the last factor
identical to that in \eqref{eq:mir_F0}. Removing that factor produces a
slightly smaller example, with only three singular fibers.

\begin{example} \label{ex:mirr3}
The factorizations $\tau_{-3,1}\cdot\tau_{0,1}\cdot \tau_{1,0}=
(a^{-3}ba^3)\cdot b\cdot a$ and $\tau_{-2,1}\cdot \tau_{0,1}\cdot
\tau_{0,1}=(a^{-2}ba^2)\cdot b\cdot b$ of the same element of $\Gamma$
are not related by Hurwitz and conjugation equivalence. They
describe inequivalent genus 1 Lefschetz fibrations with three
singular fibers, distinguished
by their first homology groups (0 vs.\ $\Z/2$).
\end{example}

In fact, we can again perform a Hurwitz move to pull out the common factor
$\tau_{0,1}=b$ from the factorizations in Example \ref{ex:mirr3}. This
yields a pair of factorizations consisting of just two Dehn twists, which
are not Hurwitz equivalent (but are related by global conjugation by $a$).

\begin{example} \label{ex:trade2}
For all $k\in\Z$, the element $a^{-4}(ab)^3\in \Gamma$ can be factored as
$\tau_{k,1}\cdot \tau_{k+2,1}$. These factorizations represent two distinct
Hurwitz equivalence classes depending on the parity of $k$, since a
Hurwitz move transforms $\tau_{k,1}\cdot \tau_{k+2,1}$ into 
$\tau_{k-2,1}\cdot \tau_{k,1}$. On the other hand, they are global
conjugates of each other (by powers of $a$).
\end{example}

In fact, \eqref{eq:mir_F1} can be rewritten in the form
$\tau_{-3,1}\cdot\tau_{-1,1}\cdot \tau_{0,1}\cdot \tau_{2,1}$ by Hurwitz
moves, so \eqref{eq:mir_F1} and \eqref{eq:mir_F0} both arise as 
fiber sums of two of the factorizations in Example \ref{ex:trade2}. In
particular, they are related by Hurwitz moves and a {\em partial conjugation}
(affecting two factors), and the corresponding 4-manifolds are related by
a Luttinger surgery~\cite{AurouxMCG}.

It is not hard to find other instances of pairs of factorizations describing
inequivalent genus 1 Lefschetz fibrations distinguished by their first
homology groups, as in Example \ref{ex:mirr3}. Perhaps more interesting is the
existence of elements of $\Gamma$ that can be factored in more than two
inequivalent ways, or of examples that can be
distinguished only by more subtle invariants. We now give a few such
examples (found by a computer search):

\begin{example}\label{ex:567}
The identities
$$T_{1,1}\cdot T_{8,-3}\cdot T_{7,-3}=T_{1,2}\cdot T_{3,1}\cdot T_{3,-1}=
T_{1,3}\cdot T_{2,1}\cdot T_{3,-1}=\mm{9&19\\44&93}$$
in $SL(2,\Z)$ lift to three factorizations of the same element of $\Gamma$
into products of Dehn twists which all belong to different Hurwitz and
conjugation equivalence classes. The first two correspond to Lefschetz
fibrations whose total space is simply connected, while the third has
a total space with first homology group $\Z/5$. All three are distinguished
by the invariant defined in \S \ref{s:invt3fibers} below, as they correspond
to the three different minimal triples $(11,10,3)$, $(5,7,6)$, and $(5,10,5)$
respectively (see \S \ref{s:invt3fibers}).
\end{example}

While this gives rise to three different Lefschetz fibrations filling the same
contact 3-manifold, it is not clear to us whether the two simply connected
fillings are different Stein manifolds, or the same manifold carrying two
different genus 1 Lefschetz fibrations. (Note that these fillings have the same signature $+1$ and first Chern class
$c_1=0$.) It is also
natural to ask whether the symplectic surfaces in $B^4$ determined by
these factorizations (viewed as products of half-twists in the braid group
$B_3$) are distinguished by the fundamental groups of their complements.
(This appears likely, but we have not been able to prove it.)

Note that the latter two of the factorizations in Example \ref{ex:567}
share the same third factor; looking only at the first two factors, we
have the following:

\begin{lemma}\label{ex:trade5}
The two factorizations in $\Gamma$ corresponding to the identities
$$T_{1,2}\cdot T_{3,1}=T_{1,3}\cdot T_{2,1}=\mm{1&-5\\5&-24}$$
are not related by Hurwitz and conjugation equivalence.
The total spaces of the corresponding Lefschetz fibrations over the disk 
are related by a ``complex conjugation'', i.e.\ there is a diffeomorphism 
between them which lifts an orientation-reversing diffeomorphism of the 
disk and maps fibers to fibers in an orientation-reversing manner.
\end{lemma}

\proof In both cases, we have a product of two Dehn twists
$\tau_1\cdot \tau_2$ about loops $\gamma_1,\gamma_2$ with
algebraic intersection number $\gamma_2\cdot\gamma_1=+5$. Thus,
there exists an oriented basis $(u,v)$ of $H_1(F,\Z)$ in which $[\gamma_2]=u$
and $[\gamma_1]=5v+ku$ for some $k\in \Z$; i.e., the factorization is
globally conjugate to $\tau_{k,5}\cdot \tau_{1,0}$. Since a change of basis (keeping
$[\gamma_2]=u$) modifies the integer $k$ by a multiple of 5, the
classification up to global conjugation is given by the
various possible values of $k$ mod 5.
In our case, we find that $\tau_{1,2}\cdot \tau_{3,1}$ is conjugate to 
$\tau_{2,5}\cdot \tau_{1,0}$, while
$\tau_{1,3}\cdot \tau_{2,1}$ is conjugate to $\tau_{3,5}
\cdot\tau_{1,0}$.

Observe now that an (inverse) Hurwitz move rewrites $\tau_{2,5}\cdot
\tau_{1,0}$ into $\tau_{1,0}\cdot \tau_{-3,5}$, and conjugating 
the factors of this latter expression by $\binom{-3\ \ 1}{\ 5\ -2}$ yields 
back $\tau_{2,5}\cdot \tau_{1,0}$. Hence, the Hurwitz equivalence
class of $\tau_{2,5}\cdot \tau_{1,0}$ consists entirely of factorizations
that are globally conjugate to it. By a similar argument, the same holds for $\tau_{3,5}\cdot
\tau_{1,0}$. Thus the integer $k$ mod 5 distinguishes the Hurwitz and
conjugation equivalence classes of the two factorizations under
consideration.

Finally, we observe that the orientation-reversing involution 
$C=\binom{0\ \,1}{1\ \,0}$ conjugates $T_{1,2}$ to $T_{2,1}^{-1}$ and
$T_{3,1}$ to $T_{1,3}^{-1}$. Thus, simultaneously applying $C$ to the 
fibers of the first Lefschetz fibration (with monodromy $\tau_{1,2}\cdot
\tau_{3,1}$) and reversing the orientation of its base (which turns 
the monodromies into their inverses) yields the second fibration (with
monodromy factorization $\tau_{1,3}\cdot \tau_{2,1}$). (A key
feature that makes the construction work is that the global monodromy
$\binom{1\ \ -5}{5\ -24}$ and its inverse are conjugate under $C$.)
\endproof

\begin{example}\label{ex:2511}
The three factorizations in $\Gamma$ corresponding to the identities
$$T_{1,1}\cdot T_{2,-3}\cdot T_{3,1}=T_{2,5}\cdot T_{1,0}\cdot T_{3,1}
=T_{3,8}\cdot T_{0,1}\cdot T_{2,1}=\mm{23&-101\\64&-281}$$
in $SL(2,\Z)$ belong to three different Hurwitz and
conjugation equivalence classes. The total spaces of the corresponding
Lefschetz fibrations are all simply connected; the Lefschetz fibrations are
distinguished by the corresponding minimal (or small) triples, which are respectively
$(5,2,-11)$, $(5,13,-1)$, and $(-3,13,2)$ (see \S \ref{s:invt3fibers}).

The first two of these factorizations have the same third factor and
differ by applying (a conjugate of) the modification described in
Lemma \ref{ex:trade5} to the first two factors. Meanwhile, the first and
third factorizations differ by the modification of Example~\ref{ex:trade2} 
(the first expression can be rewritten as
$T_{3,8}\cdot T_{1,1}\cdot T_{3,1}$ by a Hurwitz move).
\end{example}

\begin{example}\label{ex:259}
The four factorizations in $\Gamma$ corresponding to the identities
$$T_{1,3}\cdot T_{1,5}\cdot T_{2,1}=T_{2,7}\cdot T_{-1,1}\cdot T_{1,1}
=T_{2,7}\cdot T_{0,1}\cdot T_{2,1}=T_{1,2}\cdot T_{-2,1}\cdot T_{3,1}=
\mm{13&-56\\49&-211}$$
in $SL(2,\Z)$ belong to four different Hurwitz and
conjugation equivalence classes. The Lefschetz fibrations corresponding
to the first two factorizations are simply connected, and related by
a complex conjugation (i.e., applying the orientation-reversing involution
$\binom{1\ \ 0}{4\ -1}$ to the fibers and reversing the orientation of the
base). The other two factorizations correspond to Stein fillings whose
first homology groups are $\Z/2$ and $\Z/5$ respectively. The four
Lefschetz fibrations are distinguished by their minimal triples,
which are respectively $(-2,5,9)$, $(-9,5,2)$, $(-2,12,2)$ and $(-5,5,5)$.

(Note: the first and fourth factorizations are related by applying Hurwitz
moves and the modification of Lemma \ref{ex:trade5}; the second and fourth
ones as well; whereas the second and third are related by the modification
of Example \ref{ex:trade2}.)
\end{example}

\section{An invariant of Lefschetz fibrations with three singular fibers}
\label{s:invt3fibers}

Let $f:M\to D^2$ be a Lefschetz fibration with three singular fibers.
Choose a reference fiber $F$ and a distinguished collection of vanishing
paths, which determines a basis of vanishing cycles
$(\gamma_1,\gamma_2,\gamma_3)$ in $F$. Also pick arbitrary orientations of
the vanishing cycles. Then we can associate to this data the triple of
algebraic intersection numbers
$(x,y,z)=(\gamma_2\cdot\gamma_1,\gamma_3\cdot\gamma_1,\gamma_3\cdot\gamma_2)
\in\Z^3$.  We now study the dependence of this triple on the choices made,
and define an equivalence relation on $\Z^3$ so that the equivalence class
of $(x,y,z)$ is an invariant of $f$.

Changing the choices of orientations of one of the vanishing cycles changes
the signs of two of the elements in the triple, i.e.\ we have
\begin{equation}\label{eq:signs}
(x,y,z)\sim (-x,-y,z)\sim (-x,y,-z)\sim (x,-y,-z).
\end{equation}
More important is the effect of changing the distinguished basis of
vanishing paths (i.e., performing Hurwitz moves on the monodromy
factorizations). 

Since
$[\tau_{\gamma_1}(\gamma_2)]=[\gamma_2]+(\gamma_1\cdot\gamma_2)[\gamma_1]$,
replacing $(\gamma_1,\gamma_2,\gamma_3)$ by
$(\tau_{\gamma_1}(\gamma_2),\gamma_1,\gamma_3)$ changes $(x,y,z)$ to
$(-x,z-xy,y)$, or changing the orientations,
\begin{equation}\label{eq:sigma1}
(x,y,z)\sim (x,xy-z,y).
\end{equation}
Similarly, replacing $(\gamma_1,\gamma_2,\gamma_3)$ by $(\gamma_1,
\tau_{\gamma_2}(\gamma_3),\gamma_2)$ yields the triple $(y-zx,x,-z)$, or
after an orientation change,
\begin{equation}\label{eq:sigma2}
(x,y,z)\sim (zx-y,x,z).
\end{equation}
Performing both in sequence, we get
$(x,y,z)\sim (x,xy-z,y)\sim (xy-(xy-z),x,y)=(z,x,y)$.
Thus \begin{equation}\label{eq:cycle} (x,y,z)\sim (z,x,y),\end{equation}
i.e.\ triples related by cyclic permutation are equivalent.
Using cyclic permutations, we can rewrite \eqref{eq:sigma1} and
\eqref{eq:sigma2} in a slightly more symmetric manner, to yield operations
that we call {\em mutations} (in the first, second, or third position of the
triple):
\begin{equation}\label{eq:mutate}
(x,y,z)\sim (yz-x,z,y)\sim (z,xz-y,x)\sim (y,x,xy-z).
\end{equation}
Note that these operations, which modify one element of the triple while
switching the two others, are involutive. To summarize:
\begin{proposition}
The equivalence class of the triple $(x,y,z)$ up to sign changes
\eqref{eq:signs}, cyclic permutations \eqref{eq:cycle}, and mutations
\eqref{eq:mutate} is an invariant of $f$.
\end{proposition}

We now seek to find a preferred representative of a given equivalence class.

\begin{definition}
A triple $(x,y,z)$ is {\em small} if one of its elements is $-1$, $0$ or
$1$. A triple with $|x|,|y|,|z|\ge 2$ is {\em minimal (}resp.\ {\em weakly
minimal)} if \/
$|yz-x|>|x|$, $|xz-y|>|y|$, and $|xy-z|>|z|$ $($resp.\ $|yz-x|\ge|x|$, etc.$)$
\end{definition}

\noindent
In fact, a non-small triple is minimal if and only if either $xyz<0$ or
the largest of $|x|,|y|,|z|$ is less than half of the product of the other two.

\begin{proposition}
A given equivalence class either contains small representatives, or it
contains exactly one weakly minimal triple up to sign changes and
permutations (cyclic permutations only if the weakly minimal triple is
minimal, all permutations otherwise).
These possibilities are mutually exclusive.
\end{proposition}

\proof

Assume $(x,y,z)$ is a minimal triple:
then mutating in any of the three positions replaces one of $x,y,z$ 
by a new element that is the largest of the new triple, and larger than
half of the product of the two other elements. Indeed, if say we replace
$z$ by $\hat{z}=xy-z$, minimality implies that $|\hat{z}|>|z|$, which in
turn implies that $|\hat{z}|>|xy|/2$, and in particular $|\hat{z}|>|x|,|y|$.
The new triple is neither small nor weakly minimal, since mutating again
in the same place yields back the original triple.

Consider now a triple $(x,y,z)$ that is neither small nor weakly minimal, 
with say $|z|>|xy|/2$ the largest element (for instance a triple
obtained by mutating a minimal triple in the third position). Then $|x|<|yz|/2$ and $|y|<|xz|/2$. 
Thus, mutating in the first (resp.\ second) position replaces $x$ (resp.\
$y$) by $\hat{x}=yz-x$ (resp.\ $\hat{y}=xz-y$), which is the largest element 
of the new triple, as $|\hat{x}|>|yz|/2$ (resp.\
$|\hat{y}|>|xz|/2$).  However, mutating in the third position
causes that element to decrease.

Thus, if we start from a minimal triple and perform successive mutations
without ever backtracking (mutating twice the same position), we keep
obtaining larger and larger triples that are not weakly minimal, and in which the element
last modified is the largest (and larger than half of the product of the two 
others). 
This ensures that only one minimal triple (up to
cyclic permutations and sign changes) exists in the equivalence class, and
no small triples are encountered.

If the initial triple is weakly minimal but not minimal, the argument
proceeds similarly, except one of the mutations
leads to an equality. If say $|xy-z|=|z|$, then we
must have $z=xy/2$ (since $xy\neq 0$), and the mutation takes
$(x,y,z)$ to $(y,x,z)$. All other mutations lead to triples that are not
weakly minimal, with the newly modified element the largest of the triple. 
Thus, arguing as in the minimal case, the only
weakly minimal elements in the equivalence class are permutations (not
necessarily cyclic) and sign changes of $(x,y,z)$, and there are no small
elements.

Finally, given any initial triple, if it is neither small nor weakly
minimal then some mutation replaces it by a smaller triple in the same
equivalence class, and repeating the process we eventually find either a
small triple or a weakly minimal one.
\endproof

\begin{corollary}
Two Lefschetz fibrations with three singular fibers which correspond to
different minimal triples (not related by sign changes and cyclic
permutations) are not isomorphic.
\end{corollary}

Note that the monodromy factorizations in Examples \ref{ex:567},
\ref{ex:2511} and \ref{ex:259} are already given in a form where the
corresponding triples are minimal (or small, in the case of 
$(5,13,-1)$ in Example \ref{ex:2511}).

Two comments are in order. First, while this invariant can be defined
regardless of the genus of the fiber $F$, it is clearly a lot more powerful
in the genus 1 case, where the vanishing cycles are determined by their
intersection numbers up to a finite ambiguity. Second, this is an invariant
of Lefschetz fibrations but not necessarily of their total spaces, i.e.\ it
is not obvious that the Stein fillings corresponding to the various examples 
in \S \ref{s:example3fibers} are all pairwise different.

\section{Stabilization by handle attachments}

We now explain how the results in \cite{stabslf} imply the
following statement, according to which the examples discussed
in the preceding sections are intrinsically ``unstable''.

\begin{theorem}\label{thm:stab}
Let $M_1,M_2$ be two Stein fillings of the same contact 3-manifold $N$,
with the same Euler characteristic and signature.
Then there exists an exact cobordism $W$ between $N$ and some other contact
manifold $N'$ (consisting only of standard Weinstein handles) such that
attaching $W$ to $M_1$ and $M_2$ yields deformation equivalent Stein
fillings of $N'$: $M_1\cup_\partial W\simeq M_2\cup_\partial W$.
\end{theorem}

\proof
By the work of Loi-Piergallini \cite{LP} and Akbulut-Ozbagci \cite{AkOz},
the Stein fillings $M_1$ and $M_2$ carry Lefschetz fibrations $f_1,f_2$ over the disk whose
boundary is an open book compatible with the contact structure on $N$.
By a result of Giroux (Theorem~4 of~\cite{Giroux}), the open books induced by $f_1$ and
$f_2$ on $N$ have a common positive stabilization. 

Hence, after repeatedly
stabilizing $f_i$, i.e.\ attaching a 1-handle to the fiber and adding a 
new singular fiber whose vanishing cycle runs once through the new handle
(which preserves the total space $M_i$, since the new $1$- and $2$-handles
form a cancelling pair), we can ensure that the fibers of $f_1$ and $f_2$ are diffeomorphic, and
the open books induced by $f_1$ and $f_2$ on 
$N=\partial M_1=\partial M_2$ are isotopic. 

Stabilizing further if needed, we can also ensure that the fibers of $f_1$
and $f_2$ have connected boundary. Thus, the monodromies of $f_1$ and $f_2$
are described by factorizations $\F_1,\F_2$ of the same element $\phi$ as products of Dehn
twists in the mapping class group $\Map_{g,1}$ of a genus $g$ surface with
one boundary component. Moreover, the vanishing cycles of $f_1$ and
$f_2$ are all non-separating, by exactness of $M_1$ and $M_2$.

Theorem 10 of \cite{stabslf} then implies the existence of integers
$a,b,c,d,k,l$ and standard factorizations 
$\mathcal{A},\mathcal{B},\mathcal{C},\mathcal{D}$ in $\Map_{g,1}$
such that the factorizations $\F_1\cdot (\mathcal{A})^a\cdot (\mathcal{B})^b\cdot
(\mathcal{C})^c\cdot (\mathcal{D})^d$ and
$\F_2\cdot (\mathcal{A})^{a+l}\cdot (\mathcal{B})^{b-l}\cdot         
(\mathcal{C})^{c+k}\cdot (\mathcal{D})^{d-k}$
are Hurwitz equivalent. 

The Lefschetz fibrations $\hat{f}_1:\hat{M}_1\to D^2$ and
$\hat{f}_2:\hat{M}_2\to D^2$ represented by these factorizations are
isomorphic, and hence correspond
to two deformation equivalent Stein fillings $\hat{M}_1\simeq \hat{M}_2$ 
of a new contact manifold $N'$, 
obtained by attaching Weinstein handles to $f_1$ and $f_2$. 

We now argue as in \S 5 of \cite{stabslf} to prove that $k=l=0$, i.e.\ the standard pieces attached 
to $f_1$ and $f_2$ are in fact the same. For this, we calculate and compare
the Euler characteristics and signatures of $\hat{M}_1$ and $\hat{M}_2$.
The key point is that the Lefschetz fibration corresponding to $\mathcal{A}$
has 10 more singular fibers than that corresponding to $\mathcal{B}$, hence
the Euler characteristic of its total space is higher by 10, while its
signature is lower by 6; whereas for $\mathcal{C}$ and $\mathcal{D}$ the
Euler characteristics differ by 9 and the signatures by~5 (cf.\
\cite{stabslf}). Using additivity
(or Wall's non-additivity for signature, depending on how one thinks
about the cobordism between $M_i$ and $\hat{M}_i$), we conclude that
\begin{eqnarray*}
\chi(\hat{M}_2)-\chi(\hat{M}_1)&=&\chi(M_2)-\chi(M_1)+10l-9k
\quad\text{and}\\
\sigma(\hat{M}_2)-\sigma(\hat{M}_1)&=&\sigma(M_2)-\sigma(M_1)-6l+5k
\end{eqnarray*}
(cf.\ Lemmas 15 and 16 in \cite{stabslf}).
Since $M_1$ and $M_2$ have the same signature and Euler characteristic by
assumption, and $\hat{M}_1\simeq \hat{M}_2$, we conclude that $10l-9k=0$ and
$-6l+5k=0$, hence $k=l=0$, and $\hat{M}_1$ and $\hat{M}_2$ are indeed
obtained from $M_1$ and $M_2$ by attaching the same sequence of Weinstein
handles.
\endproof

Note that, while the arguments in \cite{stabslf} can be made algorithmic and
one could determine explicit values of $a,b,c,d$ for a given pair of Lefschetz
fibrations, the construction given there is far from optimal -- as
evidenced e.g.\ by the example in \S \ref{s:exampleMS}, where a single
handle attachment suffices to make the two fillings deformation equivalent.


\begin{thebibliography}{99}

\bibitem{AkOz}
S.\ Akbulut, B.\ Ozbagci, {\sl Lefschetz fibrations on compact Stein
surfaces}, Geom.\ Topol.\ {\bf 5} (2001), 319--334; Erratum: Geom.\
Topol.\ {\bf 5} (2001), 939--945.

\bibitem{AEMS}
A.\ Akhmedov, J.\ Etnyre, T.\ Mark, I.\ Smith,
{\sl A note on Stein fillings of contact manifolds},
Math.\ Res.\ Lett.\ {\bf 15} (2008), 1127--1132.

\bibitem{AurouxMCG}
D.\ Auroux,
{\sl Mapping class group factorizations and symplectic 4-manifolds: some
open problems},
Problems on Mapping Class Groups and Related Topics,
Amer.\ Math.\ Soc., Proc.\ Symp.\ Pure Math., {\bf 74} (2006), 123--132.

\bibitem{stabslf}
D.\ Auroux,
{\sl A stable classification of Lefschetz fibrations},
Geom.\ Topol.\ {\bf 9} (2005), 203--217.

\bibitem{AKO}
D.\ Auroux, L.\ Katzarkov, D.\ Orlov,
{\sl Mirror symmetry for weighted projective planes and their noncommutative
deformations}, Ann.\ Math.\ {\bf 167} (2008), 867--943.

\bibitem{AKO2}
D.\ Auroux, L.\ Katzarkov, D.\ Orlov,
{\sl Mirror symmetry for Del Pezzo surfaces: Vanishing cycles and coherent
sheaves}, Inventiones Math. {\bf 166} (2006), 537--582.

\bibitem{AKS}
D.\ Auroux, V.\,S.\ Kulikov, V.\,V.\ Shevchishin,
{\sl Regular homotopy of Hurwitz curves},
Izv. Math. {\bf 68} (2004), 521--542.

\bibitem{Beauville}
A.\ Beauville,
{\sl Les familles stables de courbes elliptiques sur $\mathbb{P}^1$ admettant quatre
fibres singuli\`eres}, C.\ R.\ Acad.\ Sci.\ Paris S\'er.\ I Math.\ {\bf 294}
(1982), 657--660.

\bibitem{Geng}
A.\ Geng, {\sl Two surfaces in $D^4$ bounded by the same knot},
J.\ Symplectic Geom.\ {\bf 9} (2011), 119--122.

\bibitem{Giroux}
E.\ Giroux, {\sl G\'eom\'etrie de contact: de la dimension trois
vers les dimensions sup\'erieures}, Proc.\ International Congress
of Mathematicians (Beijing, 2002), Vol.\ II, Higher Ed.\ Press, 2002, pp.~405--414.

\bibitem{KK}
V.\ Kharlamov, V.\ Kulikov, {\sl On braid monodromy
factorizations}, Izv.\ Math.\ {\bf 67} (2003), 79--118.

\bibitem{LP}
A. Loi, R. Piergallini, {\sl Compact Stein surfaces with boundary as
branched covers of $B^4$}, Invent.\ Math.\ {\bf 143} (2001), 325--348.

\bibitem{M} R.\ Miranda, {\sl The basic theory of elliptic surfaces},
Dottorato di Ricerca in Matematica, ETS Editrice, Pisa, 1989
({\tt http://www.math.colostate.edu/\~{}miranda/BTES-Miranda.pdf})

\bibitem{MP}
R.\ Miranda, U.\ Persson, {\sl On extremal rational elliptic surfaces},
Math.\ Z.\ {\bf 193} (1986), 537--558.

\bibitem{MP2}
R.\ Miranda, U.\ Persson, {\sl Configurations of $I_n$ fibers on
elliptic K3 surfaces}, Math.\ Z.\ {\bf 201} (1989), 339--361.

\bibitem{Moishezon}
B.\ Moishezon, {\sl Complex surfaces and connected sums of complex projective
planes}, Lecture Notes in Math.\ {\bf 603}, Springer, Heidelberg, 1977.

\bibitem{Orevkov}
S.\ Yu.\ Orevkov, {\sl Realizability of a braid monodromy by an algebraic
function in a disk}, C.\ R.\ Acad.\ Sci.\ Paris S\'er.\ I Math.\ {\bf 326}
(1998), 867--871.

\bibitem{OS03}
B.\ Ozbagci, A.\ Stipsicz, {\sl Contact 3-manifolds with infinitely many
Stein fillings}, Proc.\ Amer.\ Math.\ Soc.\ {\bf 132} (2004), 1549--1558.

\bibitem{SeBook}
P.\ Seidel, {\it Fukaya categories and Picard-Lefschetz theory},
Zurich Lect.\ in Adv.\ Math., European Math.\ Soc., Z\"urich, 2008.

\bibitem{Cadavid}
J.\,D.\ V\'elez, L.\,F.\ Moreno, C.\ Cadavid,
{\sl Hurwitz complete sets of factorizations in the modular group and the
classification of Lefschetz elliptic fibrations over the disk},
arXiv:1309.5871.

\end{thebibliography}
\end{document}